\documentclass{amsart}

\usepackage{graphicx,tikz}
\usepackage{a4wide}

\theoremstyle{theorem}

\theoremstyle{definition}

\begin{document}

\title[]{The Universal Aesthetics of Mathematics}
\markright{Mathematical synaesthesia}
\author[]{Samuel G. B. Johnson}
\address{Division of Marketing, Business \& Society, School of Management, University of Bath, UK}
\email{sgbjohnson@gmail.com }

\begin{abstract} 
The unique and beautiful character of certain mathematical results and proofs is often considered one of the most gratifying aspects of engaging with mathematics.
 We study whether this perception of mathematical arguments having an intrinsic 'character' is subjective or universal -- this was done by having test subjects with varying degrees of mathematical experience match mathematical arguments with
paintings and music: 'does this proof \textit{feel} more like Bach or Schubert?' The results suggest that such a universal 
connection indeed exists.
\end{abstract}

\author[]{Stefan Steinerberger}
\address{Department of Mathematics, Yale University, New Haven, USA}
\email{stefan.steinerberger@yale.edu}

\maketitle

\noindent

\section{Introduction}

It is not a surprising claim that the search for beauty, both in Theorems and in Proofs, is one of the great pleasures of engaging with mathematics. Indeed, many mathematicians have remarked on this connection and often the similarity to beauty in the visual arts or music is made explicit.

\begin{quote}
The mathematician's patterns, like those of the painter's or the poet's, must be beautiful; the ideas, like the colours or the words, must fit together in a harmonious way. (G.H. Hardy \cite{hardy})
\end{quote}
\begin{quote}
Why are numbers beautiful? It's like asking why is Beethoven's Ninth Symphony beautiful. If you don't see why, someone can't tell you. I know numbers are beautiful. If they aren't beautiful, nothing is. (Paul Erd\H{o}s \cite{hoffman})
\end{quote}
\begin{quote}
A scientist worthy of the name, above all a mathematician, experiences in his work the same impression as an artist; his pleasure is as great and of the same nature. (H. Poincar\'{e} \cite{poincare})
\end{quote}

Theorems can be 'deep', 'profound', 'surprising' or 'derivative' and 'boring'; conjectures can be 'daring', 'bold', 'natural' and sometimes 'false for trivial reasons' \cite{grothendieck}. Proofs can be 'beautiful', 'unexpected', 'clean', 'technical', 'elementary', 'lovely', 'nifty', 'hand-wavy' or even 'impudent' (Littlewood's description \cite{littlewood} of Thorin's proof of the Riesz-Thorin theorem). While mathematical tastes are diverse, people within the same area tend to have form sone consensus whether an argument, theorem or conjecture is beautiful, surprising, deep or insightful.
We were interested in whether or not there was any \textit{objective} basis to this phenomenon. More specifically, do proofs really have an intrinsic 'character' that is similarly perceived by different people or is this taste a consequence of mathematical socialization? Is Thorin's proof impudent or does one learn to call it that as a part of one's education?

\section{Design of the Experiment}
\subsection{Setup} The two main challenges are as follows:
\begin{enumerate}
\item finding a way to make this effect, if it exists, quantitatively measurable 
\item and ensuring that the effect is authentic and not an effect of mathematical socialization.
\end{enumerate}
The second point rules out a great number of obvious approaches (for example, having students give descriptions of the character mathematical arguments). We chose to use a comparative
approach: participants of the study were shown four mathematical arguments
(given below) and then either shown four paintings or four pieces of music; on a scale of 0--10, we asked them to rate the similarity of the piece of mathematical reasoning and the piece
of art. 

\subsection{Beautiful reasoning.} We selected four classical pieces of elementary mathematical reasoning that were chosen because of their beautiful or surprising character and
their immediate accessibility to people without mathematical training. The four arguments will be familiar to most readers and were displayed as follows.\\

\textbf{1.} 
$$ \frac{1}{2} + \frac{1}{4} + \frac{1}{8} + \frac{1}{16} + \frac{1}{32} + \dots = 1.$$
We can see this by cutting a square with total area 1 into little pieces. 
\vspace{0pt}
\begin{figure}[h!]
\begin{tikzpicture}[scale=3.2]
\draw [ultra thick] (0,0) -- (1,0);
\draw [ultra thick] (1,1) -- (1,0);
\draw [ultra thick] (1,1) -- (0,1);
\draw [ultra thick] (0,0) -- (0,1);
\node at (0.25, 0.5) {\Huge $\frac{1}{2}$};
\draw [ultra thick] (0.5,0) -- (0.5,1);
\node at (0.75, 0.75) {\Huge $\frac{1}{4}$};
\draw [ultra thick] (0.5,0.5) -- (1,0.5);
\node at (0.63, 0.25) {\Huge $\frac{1}{8}$};
\draw [ultra thick] (0.75,0.5) -- (0.75,0);
\draw [ultra thick] (0.75,0.25) -- (1,0.25);
\node at (0.87, 0.375) { $\frac{1}{16}$};
\node at (0.87, 0.125) {\dots};
\end{tikzpicture}
\end{figure}

\textbf{2.}  A quick way of computing
$$ 1 + 2 + 3 + 4 + \dots + 98 + 99 + 100 = 5050$$
is as follows: write the total sum twice and add the columns
\begin{align*}
1 &+ ~2 + ~~3 +~~~ 4 + \dots + ~~98 + 99 + 100 \\
100 &+ 99 + 98 + 98 + \dots +~~ 3 +~~ 2 +~~ 1\\
\vspace{-5pt}
\hspace{-135pt} &\rule{5.2cm}{0.4pt}\\
101 & ~~~\hspace{7pt} 101  \hspace{7pt} 101 \hspace{7pt} 101\hspace{7pt} \dots  \hspace{7pt} 101 \hspace{7pt} 101 \hspace{7pt} 101
\end{align*}
This yields a total of 100 times 101 (giving 10100) and half of that is exactly 5050.\\

\textbf{3.} In any group of 5 people, there are two who have the same number of friends within the group. 
We can see this as follows: suppose there exists somebody who is friends with everybody else. Then every person in the group has either 1, 2, 3 or 4 friends (because everybody has at least one friend). Since there are 5 people but only 4 numbers, one number has to appear twice.
If nobody is friends with all other people, then everybody has either 0, 1, 2 or 3 friends; again, since there are 5 people, one number has to appear twice.\\

\textbf{4.} The sum of consecutive odd numbers always adds up to a square number:
\begin{align*}
1 &= 1^2\\
1+3 &= 2^2 \\
1+3+5 &= 3^2\\
1+3+5+7 &= 4^2
\end{align*}
The reason is explained in the picture below: adding the next odd number creates a suitable layer for the next square.
\vspace{-10pt}
\begin{center}
\begin{figure}[h!]
\begin{tikzpicture}[scale=0.8]
\draw [ultra thick] (0,0) -- (4,0);
\draw [ultra thick] (4,4) -- (4,0);
\draw [ultra thick] (4,4) -- (0,4);
\draw [ultra thick] (0,0) -- (0,4);
\filldraw (0.5, 0.5) circle  (0.1cm);
\filldraw (1.5, 0.5) circle  (0.1cm);
\filldraw (2.5, 0.5) circle  (0.1cm);
\filldraw (3.5, 0.5) circle  (0.1cm);
\filldraw (0.5, 1.5) circle  (0.1cm);
\filldraw (1.5, 1.5) circle  (0.1cm);
\filldraw (2.5, 1.5) circle  (0.1cm);
\filldraw (3.5, 1.5) circle  (0.1cm);
\filldraw (0.5, 2.5) circle  (0.1cm);
\filldraw (1.5, 2.5) circle  (0.1cm);
\filldraw (2.5, 2.5) circle  (0.1cm);
\filldraw (3.5, 2.5) circle  (0.1cm);
\filldraw (0.5, 3.5) circle  (0.1cm);
\filldraw (1.5, 3.5) circle  (0.1cm);
\filldraw (2.5, 3.5) circle  (0.1cm);
\filldraw (3.5, 3.5) circle  (0.1cm);
\draw [thick] (0,1) -- (1,1);
\draw [thick] (1,0) -- (1,1);
\draw [thick] (0,2) -- (2,2);
\draw [thick] (2,0) -- (2,2);
\draw [thick] (0,3) -- (3,3);
\draw [thick] (3,0) -- (3,3);
\end{tikzpicture}
\end{figure}
\end{center}
\vspace{0pt}
\subsection{Beautiful art.}
The pieces of art were selected to be superficially similar: the four pieces of music all feature classical music for solo piano while the four paintings feature landscapes. For music, we used the first 20 seconds
of the following four works, respectively.
\begin{enumerate}
\item F. Schubert: Moment Musical No. 4, D 780 (Op. 94), played by D. Fray
\item J. S. Bach: Fugue from Toccata in E Minor (BWV 914), played by G. Gould
\item A. Diabelli, Diabelli's Waltz (LvB Op 120), played by G. Sokolov
\item D. Shostakovich, Prelude in D-flat major (Op.87 No. 15), played by A. Brendle
\end{enumerate}
The four pieces are all written for solo piano but have a very different character. Schubert's piece, among the most frequently performed works of his, is often formally compared to Bach but more
romantic in style; the Bach fugue, rhythmic and fast paced, has a an 'urgent' feel to it and is often grouped with those pieces of his work that show an Italian
influence. The Diabelli Waltz, famously dismissed by Beethoven as \textit{Schusterfleck} (roughly 'cobbler's patch', dismissive), is a simple classical waltz
(Kinderman \cite{kinderman} speaks of 'the banality of the theme', which is 'trite' and 'insufferably so when repeated'). Finally, Shostakovich's piece, part of his greater
cycle that had been inspired by the 200-year-anniversary of Bach's death, has been described as 'blatantly insincere', an 'artistic non-entity' that is 'cracking jokes' \cite{mark}.

\vspace{0pt}

\begin{center}
\begin{figure}[h!]
\begin{tikzpicture}[scale=1]
\node at (-2,0) {\includegraphics[width=0.65\textwidth]{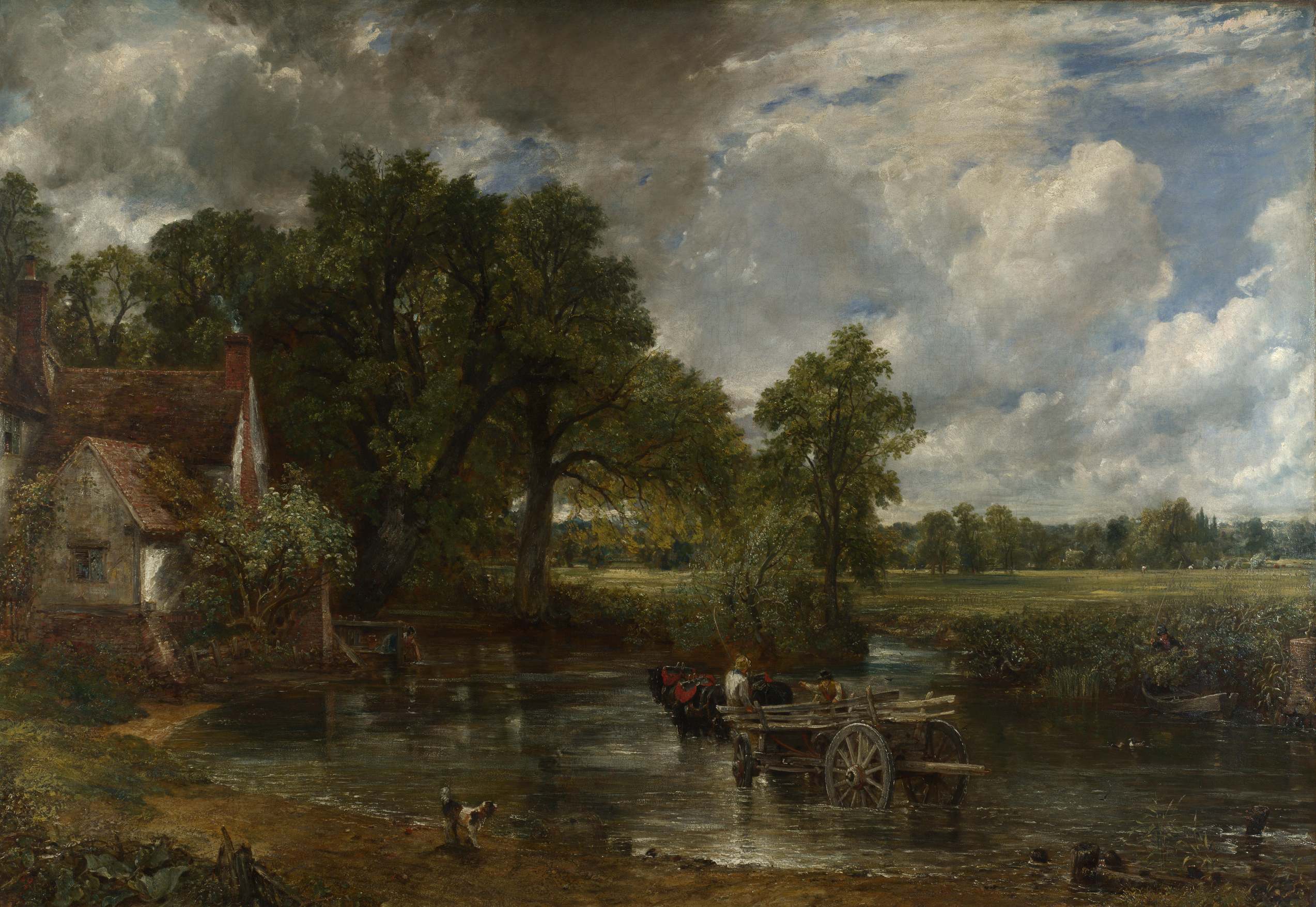}};
\end{tikzpicture}
\caption{John Constable: The Hay Wine (1821)}
\end{figure}
\end{center}
\vspace{-30pt}
\begin{center}
\begin{figure}[h!]
\begin{tikzpicture}[scale=1]
\node at (-2,0) {\includegraphics[width=0.65\textwidth]{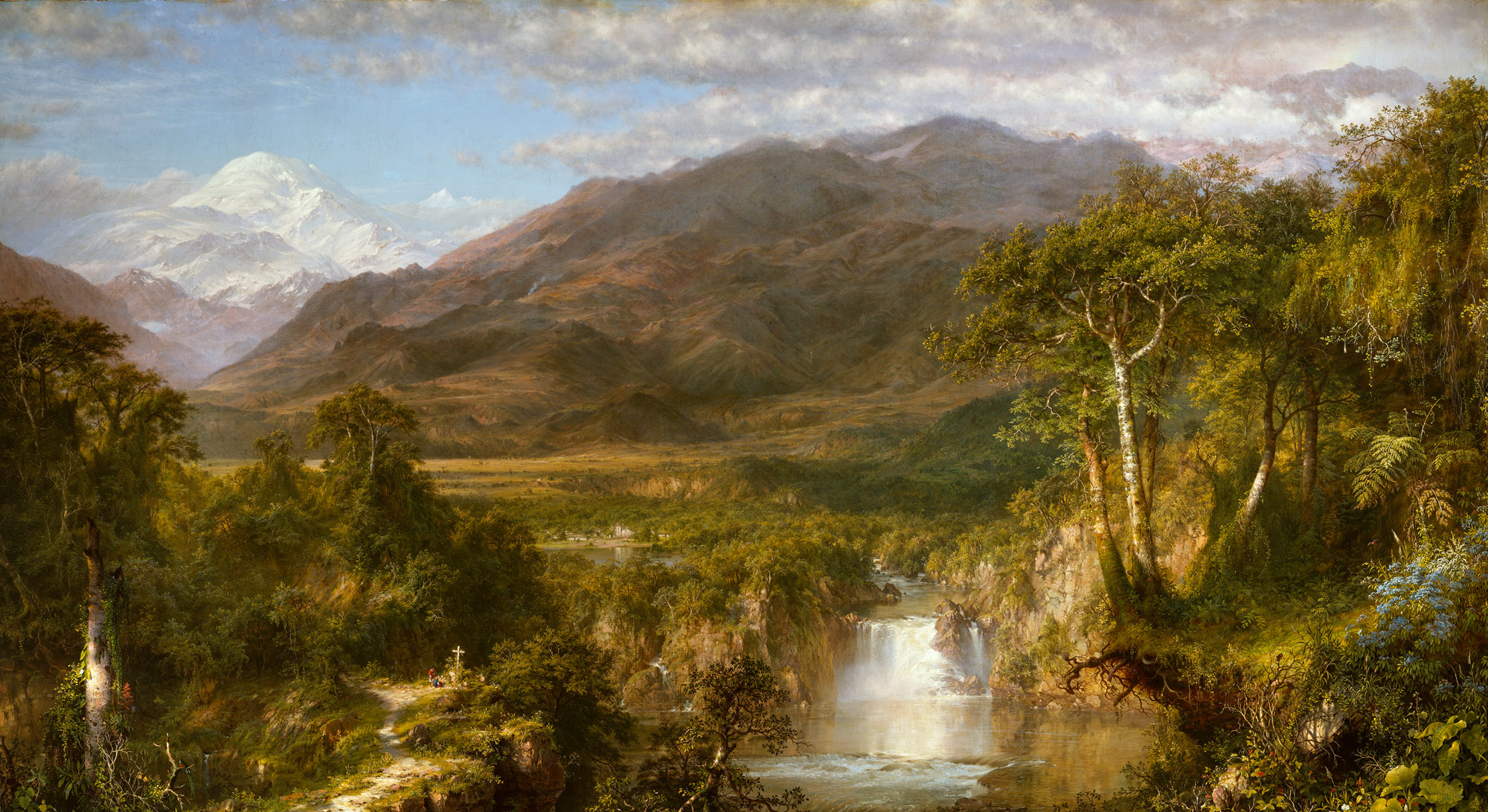}};
\end{tikzpicture}
\caption{Frederic Church: Heart of the Andes (1859)}
\end{figure}
\end{center}

\begin{center}
\begin{figure}[h!]
\begin{tikzpicture}[scale=1]
\node at (-2,0) {\includegraphics[width=0.7\textwidth]{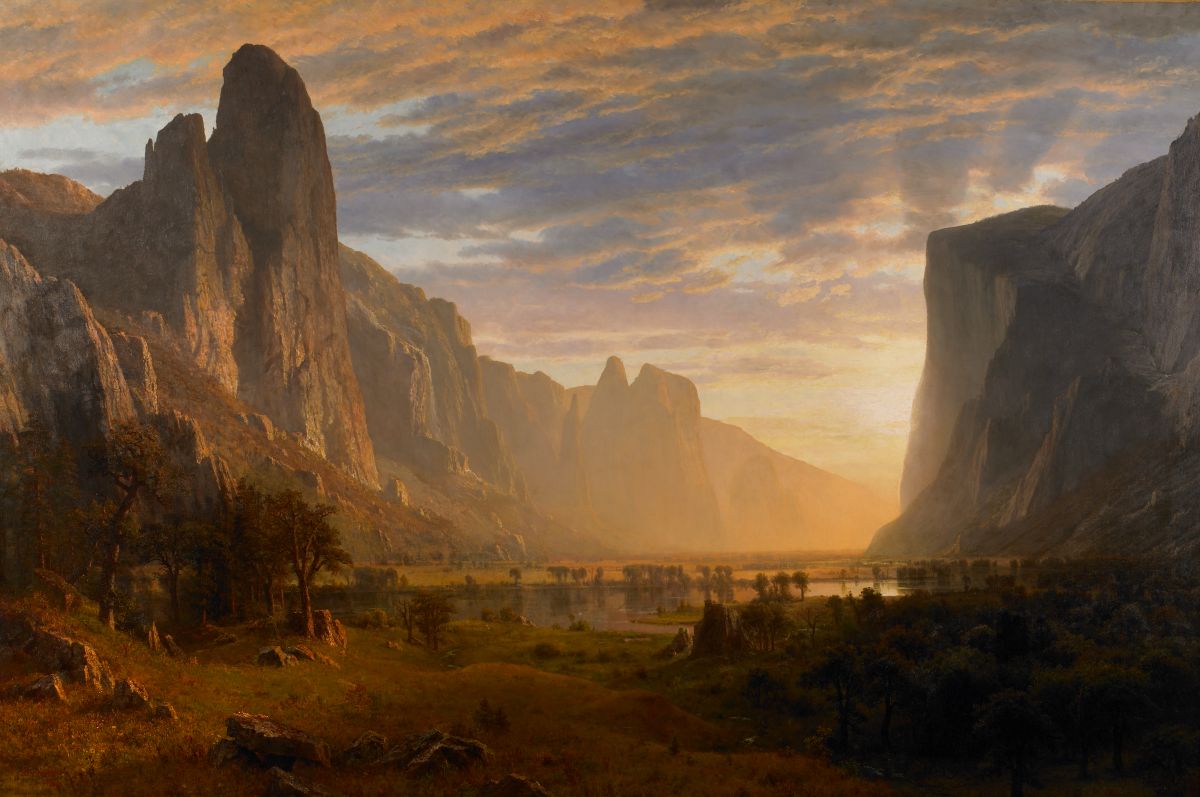}};
\end{tikzpicture}
\caption{Albert Bierstadt: Looking Down Yosemite Valley, California (1865)}
\end{figure}
\end{center}
\begin{center}
\begin{figure}[h!]
\begin{tikzpicture}[scale=1]
\node at (-2,0) {\includegraphics[width=0.7\textwidth]{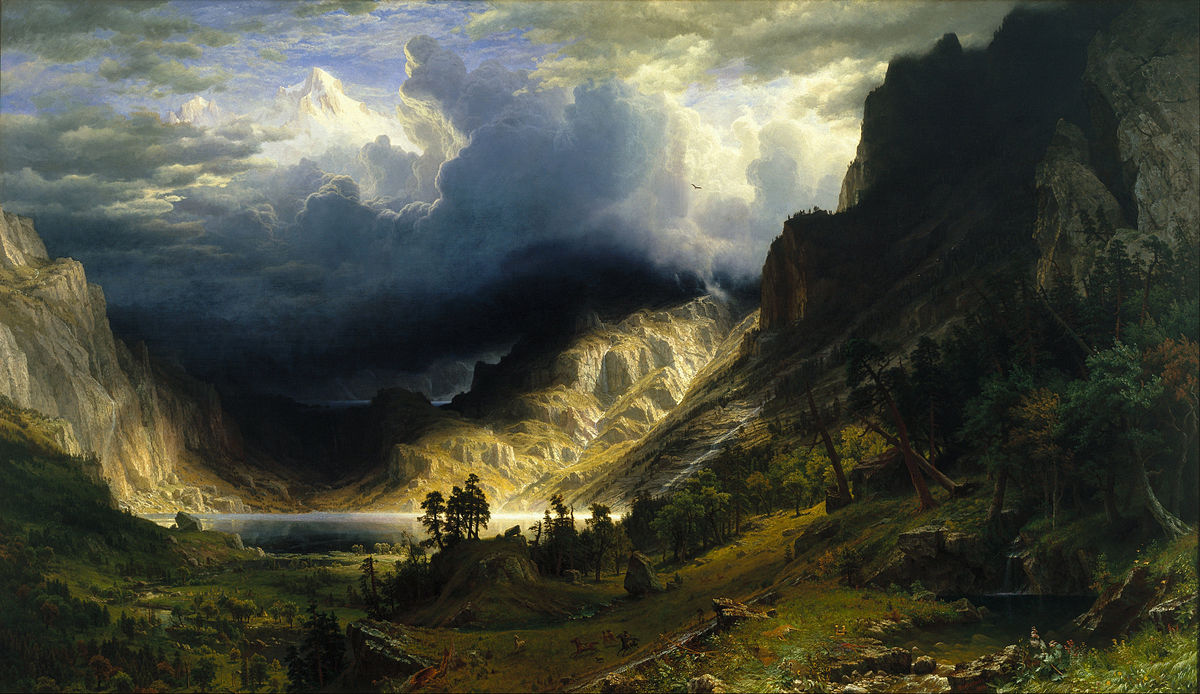}};
\end{tikzpicture}
\caption{Albert Bierstadt: Storm in the Rocky Mountains (1866)}
\end{figure}
\end{center}

We used the four paintings shown in Figures 1.-4. All of them feature realistic-romantic landscapes and date back to the 19th century.
The reader is invited to repeat the experiment: given these four pieces of music and these four paintings,
in which order do they best capture the spirit of the four mathematical arguments? In an informal survey of colleagues, most people reported feeling slightly
ill-at-ease and considered the question somewhat ill-posed.
Many of the surveyed colleagues also remarked that either matching music or matching
paintings seems easier than the other (this includes the second author who would argue that the Diabelli waltz describes
the geometric series perfectly but feels less confident in assigning a painting to it; incidentally, he is also terrible at drawing and has a hard time
appreciating the visual arts).

\section{Results.} 
\subsection{Music.} Participants were recruited from the online crowdsourcing platform Amazon Mechanical Turk ($N=299$) and were from the US; of these
$N=90$ had taken university-level math courses above Calculus. Participants read the four mathematical arguments, were then asked to reflect on the argument and
to rate the similarity to the four 20-seconds clips of music described above on a scale from 0 ('Not at all similar') to 10 ('Very similar'). The arguments as well
as the musical clips where presented in a fully randomized order and on seperate pages. After the main task, a memory check question made sure that participants
paid careful attention to the material and $N=73$ participants were dismissed at that stage. This online sample was complemented by $N=28$ Yale undergraduates and
$N=4$ professional mathematicians.
\begin{center}
\begin{figure}[h!]
\begin{tabular}{ l | c | c | c | c}
 &  Schubert & Bach & Diabelli & Shostakovich \\
Geometric Series & \textbf{4.76} &  4.39 & 4.36 & 4.62\\
Gauss Summation Trick & 4.61 &  \textbf{5.11} & 4.67 & 4.31\\
Pigeonhole Principle & 4.52 &  4.42 & \textbf{4.89} & \textbf{4.83}\\
Faulhaber Formula & 4.32 &  4.59 & \textbf{5.04} & \textbf{5.06}\\
\end{tabular}
\caption{Amazon MTurk results ($N=219$): bold denotes largest or close-to-largest degree of similarity.}
\end{figure}
\end{center}
\vspace{-10pt}
The results turned out to be far from random: taking all pairwise correlations between the $N=219$ participants of the  MTurk sample, correlation
with the mean (of all samples but the selected participant) was significantly more often positive than it was negative (145 out of 219, $p < 0.001$, $95\%-$confidence interval $(0.59, 0.72)$). This
was true for both participants who had taken higher mathematics classes $(p=0.01)$ and those who had not $(p < 0.001$). Cronbach's alpha $(\alpha =0.72)$ suggests moderate consistency across participants.
For participants with higher math training, most responses tended to correlate with both the students' mean responses (42 out of 68, $p = 0.068$) as well as the professional mathematicians' mean responses (42 out of 68, $p = 0.068$). Among participants without higher math training, however, there was no statistically reliable correlation with either the students or professionals ($p = 0.19$ and $p=0.87$, respectively).

\subsection{Paintings.} We recruited another $N=300$ participants from Amazon Mechanical Turk (of which $N=99$ had taken higher math classes and
$N=201$ had not). $N=67$ participants were dismissed after failing the memory check question and $N=1$ due to incomplete ratings. This was supplemented by $N=8$
professional mathematicians. While the connection between mathematics and music is often discussed, the same is not true for paintings and this
clearly shows in similarity ratings that were overall lower than for music. However, the association between different arguments and different works of art
is very consistent, indeed, even slightly stronger than for music. The correlation between a participant and the general mean (of all samples except the selected participant)
is positive more often than negative (156 out of 211, $p<0.001$, 95\%-confidence interval $(0.67, 0.8)$). Cronbach's alpha ($\alpha = 0.93$) indicates strongly
consistent response across participants. The agreement with professional mathematicians is also significant with 135 out of 211 correlations being positive ($p < 0.001$ and
$95\%-$confidence interval $(0.57, 0.7)$).

\begin{center}
\begin{figure}[h!]
\begin{tabular}{ l | c | c | c | c}
 &  Yosemite & Rockies & Suffolk (The Hay Wine) & Andes \\
Geometric Series & \textbf{3.51} &  2.99 & 3.30 & 3.05\\
Gauss Summation Trick & \textbf{2.38} & 2.23 & \textbf{2.43} & 1.96\\
Pigeonhole Principle & \textbf{2.42} &  2.21 & 2.25& \textbf{2.49}\\
Faulhaber Formula & 2.97 &  2.75 & \textbf{3.21} & 2.44\\
\end{tabular}
\caption{Amazon MTurk results ($N=211$): bold denotes largest or close-to-largest degree of similarity.}
\end{figure}
\end{center}

A much more extensive and complete discussion of the statistical aspects is beyond the scope of this short report and can be found in our paper \cite{paper}.

\section{Summary} We have shown that people, both with mathematical education and without, have the ability to recognize aesthetic aspects
of mathematical arguments and that these seem to be universal. The results are highly statistically significant and raise many interesting questions. 
\begin{itemize}
\item How does the effect depend on the type of art being used? People seem to see more similarity with music than with paintings but the effect seems slightly more consistent for paintings; what about sculptures, poems or possibly even jokes and puns?
\item How does the effect depend on age of the participant? How does it tie in with standard models about the development of capacity for rational thought in children?
\item How does the effect translate across cultures (where, for example, music other than the Western Canon is playing a more dominant role)?
\item Does the effect depend on the type of mathematical argument? Is an analytic inequality close to Jazz? Is a combinatorial counting argument
more of a Waltz?\\
\end{itemize}

\textbf{Acknowledgment.}
The authors wish to thank Woo-Kyoung Ahn for valuable discussions. S.S. is grateful to John Hall, Miki Havlickova,  Brett Smith and Sarah Vigliotta for helping with
the distribution of the experiment among undergraduate students.

\end{document}